\documentclass[preprint,12pt]{elsarticle}

\usepackage{amssymb}
\usepackage{amsmath}

\usepackage{color}
\usepackage{bm}
\usepackage[all]{xy}

\newtheorem{theo}{Theorem}
\newtheorem{coro}[theo]{Corollary}

\newtheorem{lemm}[theo]{Lemma}

\newenvironment{proo}[1][\proofname]{\normalfont{\itshape
		#1{:}}\quad\mdseries\ignorespaces}
{{\hspace{13cm}$\Box$}{\vskip\belowdisplayskip}}
\newcommand{\proofname}{Proof}

\newtheorem{defi}[theo]{Definition}

\journal{arXiv}

\begin{document}

\begin{frontmatter}



\title{Bi-infinite Riordan matrices: a matricial approach to multiplication and composition of Laurent series} 


%

\author[label1]{Luis Felipe Prieto-Mart\'inez} 

\affiliation[label1]{organization={Departamento de Matemática Aplicada, Universidad Politécnica de Madrid},
	addressline={Av. de Juan de Herrera 4}, 
	city={Madrid},
	postcode={28040}, 
	state={Madrid},
	country={Spain}}

\author[label2]{Javier Rico} 

\affiliation[label2]{organization={Institute for Research in Technology, ICAI School of Engineering, Comillas Pontifical University},
	addressline={C/  Rey Francisco 4}, 
	city={Madrid},
	postcode={28008}, 
	state={Madrid},
	country={Spain}}

\begin{abstract}
We propose and investigate a bi-infinite matrix approach to the multiplication and composition of formal Laurent series. We generalize the concept of Riordan matrix to this bi-infinite context, obtaining matrices that are not necessarily lower triangular and are determined, not by a pair of formal power series, but by a pair of Laurent series. We extend the First Fundamental Theorem of Riordan Matrices to this setting, as well as the Toeplitz and Lagrange subgroups, that are subgroups of the classical Riordan group. Finally, as an illustrative example, we apply our approach to derive a classical combinatorial identity that cannot be proved using the techniques related to the classical Riordan group, showing that our generalization is not fruitless.



\color{black}

\end{abstract}




\end{frontmatter}



\section{Introduction}

\subsubsection*{Basic notation}

Throughout this paper, let $\mathbb K$ be a field with identity elements denoted by $0$ and $1$. We denote by $\mathbb K[x]$ the set of polynomials with coefficients in $\mathbb K$.

\begin{itemize}
		
\item  A \textbf{formal power series} with coefficients in $\mathbb K$ is an object of the form
$$\gamma=\sum_{k=0}^\infty \gamma_kx^k.$$
		
\noindent Additionally, we denote by $\mathbb K[[x]]$ the set of formal power series with coefficients in the field $\mathbb K$ (so $\mathbb K[x]\subset\mathbb K[[x]]$).		
		
\item A \textbf{formal Laurent series} with coefficients in $\mathbb K$ is an object of the form
$$\omega=\sum_{k=-\infty}^\infty \omega_kx^k.$$

\noindent We denote by $\mathbb L(\mathbb K)$ (or simply by $\mathbb L$) the set of formal Laurent series with coefficients in $\mathbb K$ (so, again, $\mathbb K[[x]]\subset\mathbb L(\mathbb K)$).

\item We also define $\mathbb K((x))$ to be the subset of $\mathbb L(\mathbb K)$ consisting in those formal Laurent series of the form
$$\omega=\sum_{k=n}^\infty \omega_kx^k,$$
		
\noindent for some $n\in\mathbb Z$, called the \textbf{order} of $\omega$.

\item In a similar fashion, we define $\mathbb K((\frac{1}{x}))$ as the subset of $\mathbb L(\mathbb K)$ consisting in those formal Laurent series of the form
$$\omega=\sum_{k=-\infty}^n \omega_kx^k,$$

\noindent for some $n\in\mathbb Z$, that we will also  call the \textbf{order} of $\omega$.

	\end{itemize}

We use the word \emph{formal} to emphasize that the series are not assumed to be convergent, whenever this distinction is meaningful for a given choice of the field $\mathbb K$ (see the example in Section \ref{sect.comp} for a remark concerning this). From now on, we may omit the word \emph{formal}.

\subsubsection*{Multiplication and composition in $\mathbb{L}$}

The sum of elements in $\mathbb L$ is straightforward to define and is not the focus of this paper, so we omit details. However, since we will use it to define another operation below, it is worth mentioning here.

The present paper focuses on two operations in the set $\mathbb L$, namely, multiplication and composition. For that reason, both are briefly discussed next.

For two elements $\displaystyle \alpha=\sum_{i=-\infty}^\infty \alpha_{i}x^i$, $\displaystyle\beta=\sum_{j=-\infty}^\infty \beta_{j}x^j$ in $\mathbb L$, we can, under certain circumstances, define the product $\mu=\alpha\cdot \beta$ to be the Laurent series
\begin{equation} \label{eq.seriesmult}\mu = \sum_{k=-\infty}^\infty \mu_{k}x^k\text{  given, for all }k\in\mathbb Z, \text{ by } \mu_k=\sum_{i+j=k}\alpha_i\beta_j. \end{equation}

\noindent Note that, for this definition to make sense, for each $k\in\mathbb Z$, only a finite number of elements in the previous summation can be nonzero. We have the following:

\begin{lemm}[(a) is in \cite{L,N,S} and (b), (c) are just some variations] \label{lemm.mult} $ $

	\begin{itemize}
		
		\item[(a)] The set of formal power series $\mathcal F_0(\mathbb K)\subset \mathbb K[[x]]$ (in the following denoted just by $\mathcal F_0$) consisting in those elements 
		$$ \mathcal{F}_0(\mathbb{K}) = \{\omega=\omega_0+\omega_1x+\omega_2x^2+\omega_3x^3+\ldots\in \mathbb K[[x]]: \omega_0\neq 0\}$$
		
		\noindent is a group with respect to the multiplication.
		
		\item[(b)] As a consequence of the previous, $\mathbb K((x))\setminus\{0\}$ together with the multiplication is a group.
	
\item[(c)] In a similar fashion, $\mathbb K((\frac{1}{x}))\setminus\{0\}$ is also a group with respect to the multiplication (see the comment at the end of this subsection). 
	\end{itemize}

\end{lemm}

On the other hand, for two elements $\displaystyle \omega$ and $\chi$ in $\mathbb L$, we can, under certain circumstances, define the composition $\omega\circ \chi$ to be the Laurent series 
\begin{equation} \label{eq.composition}\displaystyle
	\omega \circ \chi = \sum_{i=-\infty}^\infty \omega_i \chi^i.
\end{equation}

\noindent  In general, the problem of deciding when the composition of two formal Laurent series is well defined is not immediate. For the series shown in Equation \eqref{eq.composition} to make sense, we need the power series $\ldots, \omega_{-1}\chi^{- 1}, \omega_0\chi^0,\omega_1\chi^1,\ldots$ to be well defined and the corresponding infinite sum to be well defined too. We omit details because we will discuss them later from our matricial point of view and we refer the reader to \cite{Bu, GB}.

  We have the following:

\begin{lemm}[Section 2.8 in \cite{B}] The set of formal power series $\mathcal F_1(\mathbb K)\subset \mathbb K[[x]]$ (in the following denoted just by $\mathcal{F}_1$) consisting in those elements
	$$\mathcal F_1(\mathbb K)=\{\omega=\omega_1x+\omega_2x^2+\omega_3x^3+\ldots\in \mathbb K[[x]]: \omega_1\neq 0\} $$
	
	\noindent is a group with respect to the composition. 
	
\end{lemm}

Before closing this section, it is worth to mention that there is a bijection $\mathbb K((x))\longmapsto \mathbb K((\frac{1}{x}))$ preserving the multiplication and given by 
$$\omega=\omega_nx^n+\omega_{n+1}x^{n+1}+\ldots\longmapsto \omega\circ\left(\frac{1}{x}\right)=\omega_n\frac{1}{x^n}+\omega_{n+1}\frac{1}{x^{n+1}}+\ldots$$

\subsubsection*{Infinite matrices, infinite column vectors and formal series} \label{subsection.bi}


The relation between elements in $\mathbb K[[x]]$ and lower triangular (mono-) infinite matrices has already been explored through the Riordan group. A lower triangular (mono-) infinite matrix is an object of the type
\begin{equation}\label{eq.im}[a_{ij}]_{i,j\in\mathbb N}=\begin{bmatrix}a_{0,0} \\ a_{1,0} & a_{1,1}\\ \vdots & \vdots & \ddots \end{bmatrix},\qquad a_{ij}\in\mathbb K. \end{equation}

\noindent The multiplication of these objects is well defined and the set of such matrices with non-zero elements in the main diagonal form a group with respect to this operation \cite{RS}.

For every element $(d,h)$ in $\mathcal F_0\times \mathcal F_1$ the associated \textbf{Riordan matrix} $R(d,h)$ is a lower triangular infinite matrix  $[a_{ij}]_{i,j\in\mathbb N}$, as the one in Equation \eqref{eq.im}, such that, for $j\in\mathbb N$, the column $[a_{j0},a_{j1},\ldots]^T$ contains the coefficients of $d\cdot h^j$. We have the following:

\begin{theo}[1FTRM, Theorem 3.1 in \cite{SS}] For every Riordan matrix $R(d,h)$ and for every column vector $\bm v=[v_0,v_1,\ldots]^T$, let $\bm w=R(d,h)\bm v$. If $\gamma$ is the formal power series whose coefficients are the entries in $\bm v$, then $d\cdot  (\gamma \circ h)$  is the formal power series whose coefficients are the entries in $\bm w$.
	
\end{theo}

As a consequence of the previous theorem, we can see that the set of Riordan matrices, denoted by $\mathcal R (\mathbb K)$ (or simply by $\mathcal R$), is a group with respect to the multiplication. Moreover, the multiplication of Riordan matrices is given by
$$R(d,h)R(f,g)=R(d\cdot f(h),g\circ h) $$

\noindent and the inverse of a given element is
$$(R(d,h))^{-1}=R\left(\frac{1}{d\circ h^{[-1]}},h^{[-1]}\right). $$

\noindent Both, in the previous equation as well as in the rest of the paper, we use the notation $h^{[n]}$ for the iterated composition and $h^{[-1]}$ for the compositional inverse, to distinguish from the multiplicative power and inverse, denoted by $h^n,h^{-1}$, respectively.

The Riordan group has two important subgroups that will be recalled later in this paper. The first one is the \textbf{Toeplitz subgroup}, denoted by $\mathcal T(\mathbb K)$ (or simply by $\mathcal T$), consisting in those elements in $\mathcal R$ of the type $R(d,x)$. The second one is the \textbf{Lagrange subgroup} that, for convenience, in this paper will be denoted by $\mathcal B(\mathbb K)$ (or simply by $\mathcal B$) and consists in those elements in $\mathcal R$ of the type $R(1,h)$. Note that every element $R(d,h)$ in $\mathcal R$ can be writen univocally as product of an element in $\mathcal T$ and an element in $\mathcal B$: 
$$R(d,h)=R(d,x)R(1,h).$$

In this paper, we will explore the relation between the product and composition of Laurent series and the product of (not mono-infinite matrices but) bi-infinite matrices. These matrices are objects of the type
$$[a_{ij}]_{i,j\in\mathbb Z}=\begin{bmatrix}\ddots &\vdots &\vdots &\vdots   \\ \hdots & a_{-1,-1} & a_{-1,0} & a_{-1,1} & \hdots\\ \hdots & a_{0,-1} &\fbox{$a_{0,0}$} & a_{0,1}& \hdots \\ \hdots &a_{1,-1} &a_{1,0} & a_{1,1}  & \hdots\\ & \vdots &\vdots &\vdots & \ddots \end{bmatrix},\qquad a_{ij}\in\mathbb K. $$

\noindent On the other hand, a bi-infinite column vector is an object of the type $\bm v=[\ldots, v_{-1},\fbox{$v_0$},v_1,\ldots]^T$, $v_i\in\mathbb K$.

In \cite{LMMPS} the reader may find a group of bi-infinite matrices isomorphic to the Riordan group. In the following section, we will briefly discuss some properties of the multiplication of bi-infinite matrices and the multiplication of a bi-infinite matrix with a bi-infinite column vector.

\subsubsection*{Outline of this paper}

\begin{itemize}
	
\item In Section \ref{sect.matrices} we will provide some comments and introduce notation concerning bi-infinite matrices to improve the readability of the rest of the paper.

	\item In Section \ref{sect.mult} we will study the matricial approach to the multiplication of elements in $\mathbb L$. To do so, we will assign to each element $\alpha\in\mathbb L$ a bi-infinite matrix $\bm A_\alpha$ in such a way that, when the product is defined, 
$$\bm A_\alpha\bm A_\beta=\bm A_{\alpha\cdot \beta}.$$

This leads to a first generalization of the 1FTRM (Theorem \ref{1FTRM1}) and to the description of two groups of bi-infinite matrices:
$$\{\bm A_\alpha:\alpha\in\mathbb K((x))\setminus\{0\}\},\qquad \bigg \{\bm A_\alpha:\alpha\in\mathbb K\left(\left(\frac{1}{x}\right)\right)\setminus\{0\} \bigg \}.$$

generalizing the Toeplitz subgroup of the Riordan group (Theorem \ref{toeplitz}).

	\item In Section \ref{sect.comp} we will study the matricial approach to the composition of elements in $\mathbb L$. To do so, we will assign to each element $\omega\in(\mathbb K((x))\cup\mathbb K((\frac{1}{x})))$ a bi-infinite matrix $\bm B_\omega$ in such a way that, when the product is defined, 
$$\bm B_\omega\bm B_\chi=\bm B_{\chi\circ\omega}.$$

This leads to a second generalization of the 1FTRM (Theorem \ref{1FTRM2}) and to the study of the multiplication and inversion of the  matrices $\bm B_\omega$ (Theorem \ref{lagrange}) in order to obtain a generalization of the classical Lagrange subgroup.

	\item In Section \ref{sect.riordan}, we discuss the possibility of generalizing the concept of Riordan matrix (associated to a pair of formal power series)  to the one of bi-infinite Riordan matrix (associated to a pair of formal Laurent series). To do so, we need to study when the multiplication between elements $\bm A_\alpha$ and $\bm B_\omega$ is well defined.

	\item Finally, in Section \ref{sect.app}, we present an application of the results from previous sections. The (classical) Riordan group has its origin in Combinatorics, and the generalized Riordan matrices introduced here can also be used to derive combinatorial identities.

\end{itemize}

\section{Some comments on the multiplication of bi-infinite matrices} \label{sect.matrices}

Let us start with the following:

\medskip

\noindent \textbf{Remark.} \emph{It is possible to define the product of a bi-infinite matrix $\bm M=[m_{ij}]_{i,j\in\mathbb Z}$ by an infinite column vector $\bm v=[\ldots, v_{-1},v_0,v_1,\ldots]^T$, to be the bi-infinite column vector $[\ldots,w_{-1},w_0,w_1,\ldots ]^T$ given by}
\begin{equation}  \label{eq.multvec} w_i=\sum_{j=-\infty}^\infty m_{ij}v_j, \qquad\forall i\in\mathbb Z,\end{equation}
\color{black}
\noindent  \emph{if and only if, for every $i\in\mathbb Z$, the number of nonzero products $m_{ij}v_j$ is finite. In particular, this happens if only a finite number of the elements $v_j$ are non-zero.}

\emph{In a similar fashion, for two such bi-infinite matrices  $\bm M=[m_{ij}]_{i,j\in\mathbb Z}$, $\bm N=[ n_{ij}]_{i,j\in\mathbb Z}$, the product $\bm M\bm N=[p_{ij}]_{i,j\in\mathbb Z}$,  given  by}
\begin{equation}  p_{ij}=\sum_{k=-\infty}^\infty m_{ik}n_{kj}, \qquad \label{eq.mult}\forall i,j\in\mathbb Z, \end{equation}\color{black}

\noindent  \emph{is well defined if, for every $i,j\in\mathbb Z$, the number of nonzero products $m_{ik}n_{kj}$ is finite.}

\medskip

 We have the following elementary result, that we will need later, but whose proof we will omit:

\begin{lemm} \label{lemm.bi} $ $

\begin{itemize}
\item[(a)] For bi-infinite matrices with coefficients in $\mathbb K$ we have, in some sense, an \emph{Associative Law}: if we have three bi-infinite matrices $\bm M,\bm N,\bm P$ such that the products $\bm M\bm N$, $\bm N\bm P$, $\bm M(\bm N\bm P)$, $(\bm M\bm N)\bm P$ are well defined, then $\bm M(\bm N\bm P)=(\bm M\bm N)\bm P$.

\item[(b)] In the following, $\bm I$ will denote the bi-infinite matrix $[\delta_{ij}]_{i,j\in\mathbb Z}$, where $\delta_{ij}$ stands for the Kronecker delta. For any bi-infinite matrix $\bm M$, we have $\bm M\bm I=\bm I\bm M=\bm M$.

\end{itemize}
\end{lemm}

We define \textbf{lower} (resp. \textbf{upper}) column vectors $[\ldots,v_{-1},v_0,v_1,\ldots]^T$ to be those satisfying that there exits some $k\in\mathbb Z$ such that, for  $i> k$ (resp. $i<k$), $v_i=0$.

\begin{figure}[h!]

\centering

\includegraphics[width=0.7\textwidth]{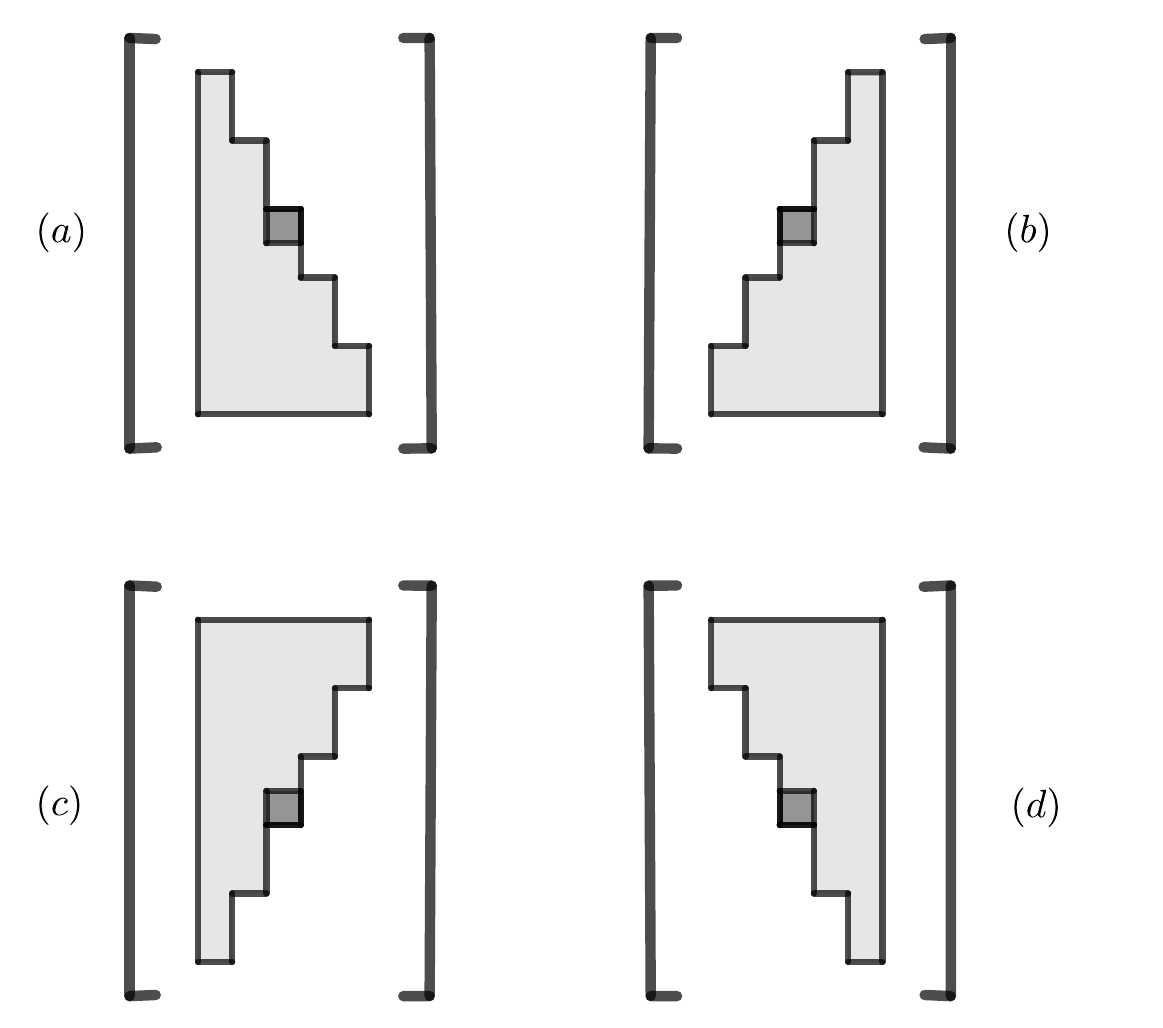}

\caption{Types of \emph{echelon forms} of bi-infinite matrices.}

\label{fig.echelon}

\end{figure}

In this paper we use some bi-infinite matrices whose shape resembles that of the \emph{reduced echelon form}. We have four cases of interest, all of them appear in Figure \ref{fig.echelon}.


We will denote by $\mathcal L_{+}$, $\mathcal L_{-}$, $\mathcal U_{-}$ and $\mathcal U_+$ to the sets of those bi-infinite matrices  corresponding to the pictures (a), (b), (c) and (d), respectively, in Figure \ref{fig.echelon}. Formally, a bi-infinite matrix $\bm M=[m_{ij}]_{i,j\in\mathbb Z}$ belongs to $\mathcal L_+$ (similar for $\mathcal L_{-}$, $\mathcal U_{-}$ and $\mathcal U_+$) if and only if  there exists a bi-infinite vector  of strictly increasing integer numbers $[\ldots,v_{-1},v_0,v_1,\ldots]$ such that 
$$m_{ij}\;\begin{cases}=0&\text{if } i<v_j\\ \neq 0 &\text{if }i=v_j\end{cases}\qquad \forall j\in\mathbb Z.$$

We have the following:

\begin{lemm} \label{lemm.matvec} Let $\bm M=[m_{ij}]_{i,j\in\mathbb Z}$ and $\bm v=[\ldots,v_{-1},v_0,v_1,\ldots]^T$.

\begin{itemize} 

\item[(a)]  If $\bm M\in \mathcal L_+$ and $\bm v$ is a lower vector, then $\bm M\bm v$ is well defined.

\item[(b)]  If $\bm M\in \mathcal L_-$ and $\bm v$ is an upper vector, then $\bm M\bm v$ is well defined.

\item[(c)]  If $\bm M\in \mathcal U_-$ and $\bm v$ is a lower vector, then $\bm M\bm v$ is well defined.

\item[(d)]   If $\bm M\in \mathcal U_+$ and $\bm v$ an upper vector, then $\bm M\bm v$ is well defined.

\end{itemize}

\end{lemm}

\begin{proo} In any of the cases above, for every $j$, the number of nonzero products $m_{ij} v_j$ in the sum in Equation \eqref{eq.multvec} is finite.

\end{proo}


Let us define the matrix $\bm J=[\xi_{ij}]_{i,j\in\mathbb Z}$ given, for every $i,j\in\mathbb Z$, by
\begin{equation} \label{eq.1/x} \xi_{ij}=\begin{cases}1 & \text{if }i+j=0\\ 0& \text{in other case} \end{cases} \end{equation}

\noindent that is, $\bm J$ is the bi-infinite matrix whose entries are 0, except for those in the secondary diagonal. This matrix will simplify our proofs from this moment on. We have the following:

\begin{lemm}   \label{lemm.geometry} $ $

\noindent (a) $\bm J$ is an \textbf{involution}, that is, $\bm J^2=\bm I$.

\noindent (b) For every bi-infinite matrix $\bm C=[c_{ij}]_{i,j\in\mathbb Z}$, the product $\bm J\bm C=[d_{ij}]_{i,j\in\mathbb Z}$ is well defined and the result is a bi-infinite matrix obtained from $\bm C$ by a \emph{reflection across the 0-row}, that is, $\forall i,j\in\mathbb Z$, $d_{ij}=c_{-i,j}$.

\noindent (c) For every bi-infinite matrix $\bm C=[c_{ij}]_{i,j\in\mathbb Z}$, the product $\bm C\bm J=[d_{ij}]_{i,j\in\mathbb Z}$ is well defined and the result is a bi-infinite matrix obtained from $\bm C$ by a \emph{reflection across the 0-column}, that is, $\forall i,j\in\mathbb Z$, $d_{ij}=c_{i,-j}$.

\end{lemm}

\begin{proo} The three statements are a consequence of the multiplication formula for bi-infinite matrices given in Equation \eqref{eq.mult}.



\end{proo}

As a consequence of the previous lemma, we have the following:

\begin{coro}\label{coro.jotas} $ $

\noindent (a) Every element in $\mathcal L_{-}$ can be writen as $\bm M\bm J$, for some $\bm M\in\mathcal L_{+}$, that is, the following sets are equal: $\mathcal L_-=\mathcal L_+\bm J$ and $\mathcal L_-\bm J=\mathcal L_+$.

\noindent (b) Every element in $\mathcal U_{+}$ can be writen as $\bm J \bm M\bm J$, for some $\bm M\in\mathcal L_{+}$, that is, the following sets are equal: $\mathcal U_+=\bm J\mathcal L_+\bm J$ and $\mathcal L_+=\bm J\mathcal U_+\bm J$.

\noindent (c) Every element in $\mathcal U_{-}$ can be writen as $\bm J\bm M$, for some $\bm M\in\mathcal L_{+}$, that is, the following sets are equal: $\mathcal U_-=\bm J\mathcal L_+$ and $\mathcal L_+=\bm J\mathcal U_-$.

This information can also be depicted in the following commutative diagram:
$$\xymatrix{L_+ \ar[r]^{R_{\bm J}} \ar[d]_{L_{\bm J}}& L_-\ar[l] \ar[d]^{L_{\bm J}}\\ U_-\ar[u] \ar[r]_{R_{\bm J}} & U_+\ar[u]\ar[l]}$$

\noindent where $L_{\bm J}$ and $R_{\bm J}$ denote the left and right multiplications by the matrix $\bm J$, respectively.

\end{coro}

Let us conclude this section with the following:

\begin{lemm} \label{lemm.tabla} Let $\bm M=[m_{ij}]_{i,j\in\mathbb Z}$ and $\bm N=[n_{ij}]_{i,j\in\mathbb Z}$ be two bi-infinite matrices. The product $\bm M\bm N$ follows the rules appearing in Table \ref{tab.lemm_table}.

\begin{table}[h]
\begin{center}
\begin{tabular}{|c| c| c| c| c|}\hline  & $\bm N\in\mathcal L_{+}$ & $\bm N\in\mathcal L_{-}$ &$\bm N\in\mathcal U_{+}$ & $\bm N\in\mathcal U_{-}$ \\ 
\hline  $\bm M\in\mathcal L_{+}$& $ (\bm M\bm N)\in\mathcal L_{+}$ & $ (\bm M\bm N)\in\mathcal L_{-}$ & not defined & not defined \\
  \hline  $\bm M\in\mathcal L_{-}$& not defined& not defined & $ (\bm M\bm N)\in\mathcal L_{-}$ & $(\bm M\bm N)\in\mathcal L_{+}$ \\ 
\hline  $\bm M\in\mathcal U_{+}$&not defined & not defined & $ (\bm M\bm N)\in\mathcal U_{+}$ & $ (\bm M\bm N)\in\mathcal U_{-}$ \\  
\hline  $\bm M\in\mathcal U_{-}$ & $ (\bm M\bm N)\in\mathcal U_{-}$ & $ (\bm M\bm N)\in\mathcal U_{+}$ & not defined& not defined\\ \hline  \end{tabular}
\end{center}
\caption{With \emph{not defined} we mean \emph{not defined in general}, noting that the intersection of the sets $\mathcal L_{+}$, $\mathcal L_{-}$, $\mathcal U_{-}$ and $\mathcal U_+$ may be non-empty.
}
\label{tab.lemm_table}
\end{table}

\end{lemm}

\begin{proo} Let us start by finding the cases in which the product $\bm M\bm N$ is well defined. To do so, we use Lemma \ref{lemm.matvec}.

Second, we need to check that, in the case $\bm M=[m_{ij}]_{i,j\in\mathbb Z},\bm N=[n_{ij}]_{i,j\in\mathbb Z}\in\mathcal L_+$, $\bm M\bm N$ belongs to $\mathcal L_+$. We will omit details, since the proof is quite similar to the one showing that the product of lower triangular matrices is lower triangular in the classical context.

Finally, to complete the information appearing in Table \ref{tab.lemm_table}, we can use Corollary \ref{coro.jotas}. Let us prove, as a matter of example, that 
$$\bm M\in \mathcal L_+,\bm N\in\mathcal L_-\implies\bm M\bm N\in\mathcal L_-$$

\noindent and we leave the rest of the cases to the reader. If $\bm N\in\mathcal L_-$, then $\bm N=\widetilde{\bm N}\bm J$, for some $\widetilde{\bm N}\in\mathcal L_+$. So $\bm M\bm N=(\bm M\widetilde{\bm N})\bm J$, where $(\bm M\widetilde{\bm N})\in \mathcal L_+$ and so $(\bm M\widetilde{\bm N})\bm J\in\mathcal L_-$.

\end{proo}

\section{First partial generalization of the 1FTRM and the Generalized Toeplitz group.} \label{sect.mult}

\begin{defi}  
	For every element $\displaystyle \alpha=\sum_{k = -\infty}^\infty \alpha_kx^k\in\mathbb L$, we define the bi-infinite matrix 
\begin{equation} \label{eq.aalpha}\bm A_{\alpha}=[a_{ij}]_{i,j\in\mathbb Z}=\left[\begin{array}{c c c c c c c}\ddots  & \vdots & \vdots& \vdots& \vdots& \vdots \\  \hdots & \alpha_0& \alpha_{-1}& \alpha_{-2}& \alpha_{-3}& \alpha_{-4}& \hdots\\
 \hdots & \alpha_1& \alpha_0& \alpha_{-1}& \alpha_{-2}& \alpha_{-3}& \hdots \\
 \hdots & \alpha_2& \alpha_1& \fbox{$\alpha_0$}& \alpha_{-1}& \alpha_{-2}& \hdots \\ 
 \hdots & \alpha_3& \alpha_2& \alpha_1& \alpha_0& \alpha_{-1}& \hdots
 \\  \hdots & \alpha_4& \alpha_3& \alpha_2& \alpha_1& \alpha_0& \hdots \\ & \vdots & \vdots& \vdots& \vdots& \vdots&\ddots \end{array}\right],\; \text{where } a_{ij}=\alpha_{i-j}. \end{equation}

\noindent  We say that this bi-infinite matrix is a \textbf{Toeplitz matrix}.
	
	
	
	
\end{defi}

Note that, if $\alpha\in\mathbb K((x))$ (resp. $\alpha\in\mathbb K((\frac{1}{x}))$), then $\bm A_\alpha\in\mathcal L_+$ (resp. $\bm A_\alpha\in\mathcal U_+$). As a consequence, we have the following:

\begin{theo}[First Partial Generalization of the 1FTRM] \label{1FTRM1} Let 
$$\alpha=\sum_{k=-\infty}^\infty \alpha_kx^k\in\mathbb L, \hspace{2cm} \beta=\sum_{k=-\infty}^\infty \beta_kx^k\in\mathbb L.$$
$$\bm v=[\ldots,\beta_{-1},\beta_0,\beta_1,\ldots]^T.$$

\begin{itemize}

\item[(a)] If $\alpha,\beta\in\mathbb K((x))$, then $\bm A_\alpha\bm v$ is well defined.

\item[(b)] If $\alpha,\beta\in\mathbb K((\frac{1}{x}))$, then $\bm A_\alpha\bm v$ is well defined.

\end{itemize}

The product is not necessarily well defined in any other case. In any of the previous cases, let $\bm w=[\ldots, w_{-1},w_0,w_1,\ldots]^T$ be the resulting vector and define the power series $\displaystyle \gamma=\sum_{k=-\infty}^\infty w_kx^k$. Then $\gamma=\alpha\cdot \beta$.

\end{theo}

\begin{proo} Statements (a) and (b) of the present theorem follow from the discussion preceding it and from Lemma \ref{lemm.matvec}.

Once this has been stated, the second part of the current proof derives from the comparison of the formula for the multiplication of two formal power series (see Equation \eqref{eq.seriesmult}) and the formula for the muliplication of a bi-infinite matrix by a bi-infinite column vector (see Equation \eqref{eq.multvec}) for the particular case of the matrix described in Equation \eqref{eq.aalpha}.




\end{proo}

As a direct consequence of the previous, we have the following

\begin{theo} \label{toeplitz} {The set of matrices $\bm A_\alpha$ such that $\alpha\in\mathbb K((x))\setminus\{0\}$ (resp. $\alpha\in\mathbb K((\frac{1}{x}))\setminus\{0\}$) is a group that will be called the \textbf{generalized lower (resp. upper) Toeplitz group}, isomorphic to the group $(\mathbb K((x)),\cdot)$. (resp. $(\mathbb K((\frac{1}{x})),\cdot)$)}.  For every $\alpha,\beta\in\mathbb K((x))$ (resp. in $\mathbb K((\frac{1}{x}))$),

\begin{itemize}

\item $\bm A_\alpha\bm A_{\beta}=\bm A_{\alpha\cdot \beta}$.
\item The identity element is $\bm A_1=\bm I$.

\item $\bm A_\alpha^{-1}=\bm A_{1/\alpha}$.

\end{itemize}
	
\end{theo}

\begin{proo} To prove the first statement, for $j\in\mathbb Z$, let us denote by $\bm v_j$ to the column vector corresponding to the $j$-th column in $\bm A_\beta$. This column vector contains the coefficients of $x^j\beta$. The $j$-th column in the product $\bm A_\alpha\bm A_{\beta}$ is $\bm A_\alpha\bm v_j$. According to  Theorem \ref{1FTRM1}, this column contains the coefficients of $x^j\alpha\beta$. This completes the proof of this statement, and the rest of them are a consequence of this it.

\end{proo}

Some further questions concerning the existence of multiplicative inverses of elements in $\mathbb L$ were already studied in \cite{Bu, GB}, but the information provided in the previous result will be sufficient for the rest of the paper.

\section{Second Partial Generalization of the 1FTRM and the Generalized Lagrange group} \label{sect.comp}

Some questions concerning the existence of the composition $\chi\circ\omega$, for $\chi,\omega\in\mathbb L$ can be found in \cite{Bu}. In any case, we will follow our own approach for this problem.

Note that $\mathbb K((x))$ is not closed under composition.

\medskip

\noindent \textbf{Example.} \emph{For instance,}
$$\chi=1+x+x^2+x^3+\ldots\in\mathbb K((x)),\; \omega=\frac{1}{x}\in\mathbb K((x)),$$
$$\chi \circ \omega=1+\frac{1}{x}+\frac{1}{x^2}+\frac{1}{x^3}+\ldots\notin \mathbb K((x)). $$
 
\bigskip

With respect to this, we also have a problem of interpretation.

\medskip

\noindent \textbf{Example (continuation).} \emph{Let us consider the previous example in the case $\mathbb K=\mathbb C$. Then, we can identify}
$$\chi=\frac{1}{1-x} $$

\noindent \emph{and then} 
$$\chi\circ \omega=\frac{1}{1-\frac{1}{x}}=-\frac{x}{1-x}=-x-x^2-x^3-\ldots$$

\noindent  \emph{So the Laurent series}
$$1+\frac{1}{x}+\frac{1}{x^2}+\frac{1}{x^3}+\ldots,\qquad -x-x^2-x^3-\ldots$$

\noindent \emph{are \emph{different as formal objects} but \emph{equal in some sense}. This is because they correspond to two different expansions of the same complex function of one complex variable. For us, these two formal Laurent series will be different and the other interpretation will be left appart in this paper.}

\medskip

\begin{defi} \label{defi.B}
	For every element $\omega\in\mathbb K((x))\cup\mathbb K((\frac{1}{x}))$, we define the matrix  $\bm B_\omega$ to be the bi-infinite matrix $\bm B_\omega=[b_{ij}]_{-\infty<i,j<\infty} $  such that, for every $\displaystyle j\in\mathbb Z$, $\displaystyle \omega^j=\sum_{i = -\infty}^{\infty}b_{ij}x^i$.

\end{defi}

The condition $\omega\in\mathbb K((x))\cup\mathbb K(\frac{1}{x}))$ ensures that the powers 
$$\ldots,\omega^{-1},\omega^0,\omega^1,\omega^2,\omega^3,\ldots$$
\noindent are well defined.

As a first example of such a matrix, note that  $\bm J=\bm B_{1/x}$ (the definition appears in Equation \eqref{eq.1/x}). On the other hand, we have that:
\begin{itemize}

\item If $\omega\in\mathbb K((x))$ and its order is positive, then $\bm B_\omega\in\mathcal L_+$.

\item If $\omega\in\mathbb K((x))$ and its order is negative, then $\bm B_\omega\in\mathcal L_-$.

\item If $\omega\in\mathbb K((\frac{1}{x}))$ and its order is positive, then $\bm B_\omega\in\mathcal U_+$.

\item If $\omega\in\mathbb K((\frac{1}{x}))$ and its order is negative, then $\bm B_\omega\in\mathcal U_-$.

\end{itemize}

We have the following:

\begin{theo}[Second Partial Generalization of the 1FTRM] \label{1FTRM2} Let 
$$\chi=\sum_{k=-\infty}^\infty \chi_kx^k\in\mathbb L(\mathbb K),\qquad\bm v=[\ldots,\chi_{-1},\chi_0,\chi_1,\ldots]^T.$$

\begin{itemize}

\item[(a)] If $\omega\in \mathbb K((x))$ and has positive order and $\chi\in\mathbb K((x))$, then  $\bm B_\omega\bm v$ is well defined.


\item[(b)] If $\omega\in \mathbb K((\frac{1}{x}))$ and has negative order and $\chi\in\mathbb K((x))$, then  $\bm B_\omega\bm v$ is well defined.



\item[(c)] If $\omega\in\mathbb K((x))$ and has negative order and $\chi\in\mathbb K((\frac{1}{x}))$,  then $\bm B_\omega\bm v$ is well defined.


\item[(d)] If $\omega\in\mathbb K((\frac{1}{x}))$ and has positive order and $\chi\in\mathbb K((\frac{1}{x}))$,  then $\bm B_\omega\bm v$ is well defined.


\item[(e)] If $\omega\in \mathbb K((x))\cup\mathbb K\left(\left(\frac{1}{x}\right)\right)$ and $\chi\in  \mathbb K((x))\cap\mathbb K\left(\left(\frac{1}{x}\right)\right)$, for some $n\in\mathbb Z$,  then $\bm B_\omega\bm v$ is well defined.

\end{itemize}

In any of the previous cases, the resulting vector $\bm w=[\ldots, w_{-1},w_0,w_1,\ldots]^T$ corresponds to the Laurent series $\psi$
$$\psi=\sum_{k=-\infty}^\infty w_kx^k=\chi\circ\omega$$
		
\noindent where this composition of elements in $\mathbb L$ is also well defined.
	
\end{theo}

\begin{proo} We can prove Statements (a), (b), (c), (d) and (e) using the discussion before this theorem and Lemma \ref{lemm.matvec}.




To prove the final statement, let $\bm B_\omega=[b_{ij}]_{i,j\in\mathbb Z}$ and let us denote by $\bm u_j$ to the column vector $[\ldots, b_{-1,j},b_{0j},b_{1j},\ldots]^T$. The column vector 
$$\bm w=\ldots+\chi_{-1}\bm u_{-1}+\chi_0\bm u_0+\chi_1\bm u_1+\ldots $$

\noindent that we have ensured that is well defined (at each level there are only finitely many non-zero entries), corresponds to the Laurent series
$$\chi\circ\omega=\sum_{k=-\infty}^\infty \chi_k\omega^k. $$

\end{proo}

Lemma \ref{lemm.geometry} has an interesting translation to this context. We have the following:

\begin{coro} For every $\omega\in\mathbb K((x))\cup\mathbb K((\frac{1}{x}))$, $\bm J\bm B_\omega=\bm B_{\omega\circ\frac{1}{x}}$, $\bm B_\omega\bm J=\bm B_{1/\omega}$.

\end{coro}

\begin{proo} This result is a direct consequence of Theorem \ref{1FTRM2}, particularly, Statement (e).

\end{proo}

\begin{theo} \label{lagrange} $ $

\begin{itemize}

\item[(a)] When defined, the product $\bm B_{\omega}\bm B_{\chi}$ equals $\bm B_{\chi\circ\omega}$.

\item[(b)] In view of the previous lemma, to decide if $\bm B_{\omega}\bm B_{\chi}$ is well defined, it suffices to know that, 

\begin{itemize} 

\item For every bi-infinite matrix $\bm I\bm C=\bm C\bm I=\bm C$.

\item $\bm J\bm J=\bm I$.

\item For every $\bm B_{\omega},\bm B_{\chi}\in\mathcal L_{+}$, the following products are well defined:
$$\bm B_\omega\bm B_\chi,\qquad \bm J \bm B_\omega\bm B_\chi, \qquad \bm B_\omega\bm B_\chi\bm J, \qquad\bm J\bm B_\omega\bm B_\chi\bm J $$

\item For every $\bm B_{\omega},\bm B_{\chi}\in\mathcal L_{+}$, in general, the product $\bm B_\omega\bm J\bm B_\chi$ is not  well defined.

\end{itemize}

Equivalently, we can use the comments after Definition \ref{defi.B} (to determine, in each case, if the corresponding matrix is in $\mathcal L_+$, $\mathcal L_-$, $\mathcal U_+$ or $\mathcal U_-$) along with Lemma \ref{lemm.tabla}.



\item[(c)] Let $\omega\in\mathbb K((x))\cup\mathbb K((\frac{1}{x}))$ not of order 0. There exists some $\chi\in\mathbb K((x))\cup\mathbb K((\frac{1}{x}))$ such that $\bm B_{\chi}$ is the inverse of $\bm B_\omega$ if and only if the order of $\omega$ is $\pm 1$. This $\chi$ is a compositional inverse of $ \omega$ in $\mathbb L$, belongs to the same set $\mathbb K((x)),\mathbb K((\frac{1}{x}))$ than $\omega$ and has the same order.

\item[(d)] The set of matrices $\bm B_\omega\in\mathcal L_{+}$ and the set of matrices $\bm B_\omega\in \mathcal U_{+}$, endowed with the multiplication, are monoids. 

\item[(e)] The set of matrices $\mathcal B_\omega\in\mathcal L_{+}$ such that $\omega$ has order 1 and the set of matrices $\mathcal B_\omega\in \mathcal U_{+}$ such that $\omega$ has order 1, endowed with the multiplication, are groups that will be called \textbf{generalized lower/upper Lagrange group}, respectively. The function given by
$$\bm B_{\omega}\mapsto \bm J\bm B_\omega\bm J,$$

\noindent is an isomorphism (from the first one to the second one and from the second one to the first one)  and an involution (the iterated composition of this isomorphism is the identity).


\end{itemize}

\end{theo}

\begin{proo} 
\noindent \fbox{Statement (a)} Let us denote by $\bm v_j$ to the $j$-th column of the matrix $\bm B_\chi$. Then $\bm B_\omega\bm v$, according to Definition \ref{defi.B} and to Theorem \ref{1FTRM2}, contains the coefficients of $(\chi\circ\omega)^j$. So the product $\bm B_{\omega}\bm B_{\chi}$ is the matrix whose $j$-column contains the coefficients of $(\chi\circ\omega)^j$ and thus, by Definition \ref{defi.B}, equals $\bm B_{\chi\circ\omega}$.

\medskip

\noindent \fbox{Statement (b)} The first statement appears in Lemma \ref{lemm.bi}. The second statement appears in Lemma \ref{lemm.geometry}. The third and fourth statements are a consequence of Theorem \ref{1FTRM2}.

\medskip

\noindent \fbox{Statement (c)} Let us prove the necessity of the condition of $\omega$ having order $\pm 1$.

\noindent Let us start with the case $\bm B_\omega\in\mathcal L_{+}$. Let us see that $\omega$ necessarily is of order $1$. Consider the bi-infinite linear system: 
\begin{equation}\label{eq.system}\bm B_\omega\begin{bmatrix}\vdots\\ \chi_{-2} \\ \chi_{-1}\\ \chi_0\\ \chi_1\\ \chi_2\\ \vdots \end{bmatrix}=\begin{bmatrix} \vdots \\ 0\\  0 \\ 0 \\ 1 \\ 0 \\ \vdots \end{bmatrix}. \end{equation}

\noindent Suppose that the order of $\omega$ is $n>1$ (we are in the Case (a) of Figure \ref{fig.echelon}). Let us consider, in turn, the following linear system (obtained deleting some rows)
\begin{equation}\label{eq.systemsimp}\begin{bmatrix}\ddots \\ \hdots & b_{-n,-1}  \\ \hdots &b_{0,-1} & 1 \\ \hdots &b_{n,-1} &0 & b_{n,1} \\ &\vdots &\vdots&\vdots&\ddots \end{bmatrix}\begin{bmatrix}\vdots \\ \chi_{-1}\\ \chi_0\\ \chi_1\\ \vdots \end{bmatrix}=\begin{bmatrix}\vdots\\ 0\\ 0 \\ 0 \\ \vdots   \end{bmatrix}. \end{equation}

\noindent Note that if the system showed in Equation \eqref{eq.system} has a solution, then the system showed in Equation \eqref{eq.systemsimp} also does. Note that, if $n=1$, the column vector is non-trivial.

\noindent  If $\chi\in\mathbb K((x))\cup\mathbb K((\frac{1}{x}))$ and has order $k$, this system reduces to one of the following types (depending on $\chi\in\mathbb K((x))$ or $\chi\in\mathbb K((\frac{1}{x}))$):
$$\begin{bmatrix} b_{kk}\\ b_{k+n,k} & b_{k+n,k+1}\\ \vdots & \vdots & \ddots \end{bmatrix} \begin{bmatrix}\chi_k\\ \chi_{k+1}\\ \vdots \end{bmatrix} =\begin{bmatrix} 0 \\ 0 \\ \vdots \end{bmatrix},\qquad \begin{bmatrix}\ddots \\ \hdots & b_{k-n,k-1} \\ \hdots & b_{k,k-1}& b_{kk}\end{bmatrix} \begin{bmatrix} \vdots \\ \chi_{k-1}\\ \chi_k \end{bmatrix} =\begin{bmatrix} \vdots \\ 0 \\ 0\end{bmatrix}.$$

\noindent Note that the matrices of coefficients in the previous equation are not bi-infinite and that the elements in the main diagonal are different from 0.

\noindent In both cases we can conclude that the only solution is $\chi=0$ (performing some kind of \emph{Gaussian elimination} that, although requires to replace a row by a sum of infinitely many rows, it is well defined), which is a contradiction, since $\bm B_\omega\bm B_0\neq \bm I$.

 For the cases $\bm B_\omega\in\mathcal L_{-}$, $\bm B_\omega\in\mathcal U_{-}$,$\bm B_\omega\in\mathcal U_{+}$, we have that
\begin{equation}\label{eq.cases}\bm B_\omega=\bm B_{\widetilde \omega}\bm J,\qquad \bm B_\omega=\bm J\bm B_{\widetilde \omega},\qquad \bm B_\omega=\bm J\bm B_{\widetilde \omega}\bm J, \end{equation}

\noindent respectively, for some $\bm B_{\widetilde \omega}\in\mathcal L_{+}$.

Let us see that this condition is sufficient for having an inverse. Suppose that $\bm B_{\omega}\in\mathcal L_{+}$ and $\omega$ has order 1. In this case, the classical theory for the group $(\mathcal F_1,\circ)$ ensures that $\omega$ has a compositional inverse $\omega^{[-1]}$ of order 1 and so $\bm B_{\omega^{[-1]}}$ is the inverse of $\bm B_\omega$. Again, using Corollary \ref{coro.jotas}, we can study the rest of the cases.

For the matrices in Equation \eqref{eq.cases}, if  $\omega$ has order $n$, then $\widetilde \omega$ has order $\pm n$. Moreover,  the matrices $\bm B_{\omega}$ in Equation \eqref{eq.cases} have an inverse if and only if so does $\bm B_{\widetilde \omega}$. This proves the necessity of the condition of $\omega$ having order $\pm 1$.








\medskip

\noindent \fbox{Statement (d)} This is a consequence of Lemma \ref{lemm.bi}.

\medskip

\noindent \fbox{Statement (e)} This is a consequence of the first part and of Statement (c). The second part is a consequence of Statement (b) and Corollary \ref{coro.jotas}.

\end{proo}

\section{Generalized Riordan Group} \label{sect.riordan}

\begin{defi} \label{defi.R} Let $\alpha,\omega$ be two elements, both in $\mathbb K((x))$ or both in $\mathbb K((\frac{1}{x}))$. We define the \textbf{(bi-infinite) Riordan matrix} matrix  $\bm R_{\alpha,\omega}$ to be the bi-infinite matrix $\bm R_{\alpha,\omega}=[r_{ij}]_{-\infty<i,j<\infty} $  such that, for every $\displaystyle j\in\mathbb Z$, $\displaystyle \alpha\cdot \omega^j=\sum_{i = -\infty}^{\infty}r_{ij}x^i$.

\end{defi}

\noindent The condition at the beginning ensures that the powers $\ldots,\frac{\alpha}{\omega^2},\frac{\alpha}{\omega},\alpha,\alpha\omega, \alpha\omega^2,\ldots$  are well defined. Note that $\bm R_{\alpha,x}=\bm A_\alpha$ and $\bm R_{1,\omega}=\bm B_\omega$.

Let us remark that the case $(\alpha,\omega)\in \mathcal F_0\times\mathcal F_1$ was already introduced in \cite{LMMPS}. In this case, the matrix $\bm R_{\alpha,\omega}$ is lower triangular and of the type
$$ \left[\begin{array}{c c c | c c c} \ddots & & & & &\\ \hdots & d_{-2,-2}& & & &\\ \hdots &d_{-1,-2}&d_{-1,-1}  & & &\\ \hline \hdots &d_{0,-2} &d_{0,-1} & & &\\  \hdots & d_{1,-2}& d_{1,-1}& & R(\alpha,\omega) &\\  & \vdots &\vdots & & &\end{array}\right]$$

\noindent where $R(\alpha,\omega)$ is a classical Riordan matrix.

\begin{lemm} \label{lemm.AB} Every Riordan matrix $\bm R_{\alpha,\omega}$ satisfies $\bm R_{\alpha,\omega}=\bm A_\alpha\bm B_\omega$. Furthermore, $\bm R_{\alpha,\omega}$ lies in any of the sets $\mathcal L_+$, $\mathcal L_-$, $\mathcal U_+$, $\mathcal U_-$ if and only if $\bm B_\omega$ does.

\end{lemm}

\begin{proo} The first part is a direct consequence of Theorem \ref{1FTRM1} and the second part ir a direct consequence of Lemma \ref{lemm.tabla}.

\end{proo}

\begin{theo}[Final Generalization of the 1FTRM] \label{1FTRM3} Let 
$$\chi=\sum_{k=-\infty}^\infty \chi_kx^k\in\mathbb L(\mathbb K),\qquad\bm v=[\ldots,\chi_{-1},\chi_0,\chi_1,\ldots]^T.$$

\begin{itemize}

\item[(a)] If $\omega\in \mathbb K((x))$ and has positive order and $\chi\in\mathbb K((x))$, then  $\bm R_{\alpha,\omega}\bm v$ is well defined.


\item[(b)] If $\omega\in \mathbb K((\frac{1}{x}))$ and has negative order and $\chi\in\mathbb K((x))$, then  $\bm R_{\alpha,\omega}\bm v$ is well defined.



\item[(c)] If $\omega\in\mathbb K((x))$ and has negative order and $\chi\in\mathbb K((\frac{1}{x}))$,  then $\bm R_{\alpha,\omega}\bm v$ is well defined.


\item[(d)] If $\omega\in\mathbb K((\frac{1}{x}))$ and has positive order and $\chi\in\mathbb K((\frac{1}{x}))$,  then $\bm R_{\alpha,\omega}\bm v$ is well defined.


\item[(e)] If $\omega\in\mathbb K((x))\cup\mathbb K \left(\left(\frac{1}{x}\right)\right)$ and $\chi\in\mathbb K((x))\cap\mathbb K \left(\left(\frac{1}{x}\right)\right)$, for some $n\in\mathbb Z$, then $\bm R_{\alpha,\omega}\bm v$ is well defined.

\end{itemize}
		
In any of the previous cases, the resulting vector $\bm w=[\ldots, w_{-1},w_0,w_1,\ldots]^T$ corresponds to the Laurent series $\psi$
$$\psi=\sum_{k=-\infty}^\infty w_kx^k=\alpha\cdot (\chi\circ\omega).$$

\end{theo}

\begin{proo} First, in view of Lemma \ref{lemm.AB}, we write $\bm R_{\alpha,\omega}=\bm A_\alpha\bm B_\omega$. Then, we use Theorem \ref{1FTRM2} to give meaning to $\bm B_\omega\bm v$. Finally, we use Theorem \ref{1FTRM1} to give meaning to $\bm A_\alpha (\bm B_\alpha\bm v)$.

\end{proo}

\begin{theo} $ $

\begin{itemize}

\item[(a)] A Riordan matrix $\bm R_{\alpha,\omega}$ is invertible if and only if $\alpha\neq 0$ and $\bm B_{\omega}$ is invertible, that is, the order of $\omega$ is $\pm 1$.

\item[(b)] The multiplication of two Riordan matrices $\bm R_{\alpha,\omega}\bm R_{\beta,\chi}$ is well defined if and only if $\bm B_\omega\bm B_\chi$ is well defined. In this case, $\bm R_{\alpha,\omega}\bm R_{\beta,\chi}=\bm R_{\alpha\cdot (\beta\circ\omega),\chi\circ\omega}$.

\item[(c)] The sets of matrices $\bm R_{\alpha,\omega}\in\mathcal L_+$ and the set of matrices $\bm B_\omega\in\mathcal U_+$ endowed with the multiplication, are monoids.

\item[(d)] The sets of matrices $\bm R_{\alpha,\omega}\in\mathcal L_+$ such that $\alpha\neq 0$ and $\omega$ has order 1 and the set of matrices $\bm R_{\alpha,\omega}\in\mathcal U_+$ such that $\alpha\neq 0$ and $\omega$ has order 1, endowed with the multiplication, are groups, that will be called the \textbf{generalized lower/upper Riordan group}.

\end{itemize}

\end{theo}

\begin{proo} First, Statement (a) is a consequence of the decomposition in Lemma \ref{lemm.AB} and of the existence of inverses appearing in Theorems \ref{toeplitz} and \ref{lagrange}.

Second, the first part of Statement (b) is a consequence of Lemma \ref{lemm.AB} and Lemma \ref{lemm.tabla}. To prove the second part, note that, according to Theorem \ref{1FTRM2}, the product $\bm B_\omega\bm A_\alpha$ is well defined and 
$\bm B_\omega\bm A_\alpha=\bm A_{\alpha\circ \omega}\bm B_\omega $.

\noindent So, 
$$\bm R_{\alpha,\omega}\bm R_{\beta,\chi}=\bm A_\alpha\bm B_\omega\bm A_\beta\bm B_\chi=\bm A_\alpha\bm A_{\beta\circ\omega}\bm B_\omega\bm B_\chi. $$

Third, Statements (c) and (d) are a combination of Theorems \ref{toeplitz} and \ref{lagrange}.

\end{proo}

\section{Application: Palyndromic Polynomials} \label{sect.app}

When the study of the Riordan group began in the paper \cite{SS}, the goal was not to analyze an algebraic abstract object, but to investigate an interesting structure due to its applications in Combinatorics. In this section, we aim to demonstrate that, like the classical Riordan group, the generalized Toeplitz, Lagrange, and Riordan groups presented here also provide a framework that facilitates the proof of various combinatorial identities, including some not available for the classical Riordan group. Furthermore, these identities can be studied within the context of the \emph{Concrete Mathematics} (in the sense of the term used by R. L. Graham, E. D. Knuth and O. Patashnik in \cite{GKP}). As an example, we present a really short proof of a well-known identity from the theory of finite simplicial complexes.

For the sake of brevity, we refer the reader to \cite{LMP} for the concepts of \emph{finite simplicial complex} and \emph{f-polynomial}. None of these concepts are crucial. We only need to know that the \emph{extended f-polynomial} of a simplicial complex is a polynomial $f(x)=f_{-1}+f_0x+\ldots+f_dx^{d+1}$, being $d$ the dimension of the simplicial complex. In the aforementioned paper \cite{LMP}, the relation of this kind of problems regarding f-vectors and Riordan matrices was explored. The reader may find more information about problems concerning f-vectors in \cite{G, Z} and the references therein.

For some families of finite simplicial complexes (like partitionable simplicial complexes), it makes sense to study, not only the extended f-polynomial but also the \textbf{h-polynomial}, defined as 
\begin{equation} \label{eq.hf} h(x)=h_0+h_1x+\ldots+h_{d+1}x^{d+1}=(1-x)^{d+1}f\left(\frac{x}{1-x}\right).\end{equation}

\noindent  (Section 8.3 in \cite{Z}), since, in this case, the coefficients $h_0,\ldots, h_{d+1}$ have combinatorial meaning.

We also need to recall that, for many important families of finite simplicial complexes (such as simplicial polytopes), the $h$-polynomial is \textbf{palyndromic}, that is, 
\begin{equation}\label{eq.palyndromic} h_k=h_{d+1-k},\qquad \forall k=0,\ldots, d+1.\end{equation}

\noindent These conditions are known as the \textbf{Dehn-Sommerville Equations} (Theorem 8.21 in \cite{Z}).

Our approach allows us to find very quickly another version of interest of the Dehn-Sommerville equations stated only in terms of the f-vector that can be found in Section 9.2 in \cite{G}.

\begin{theo} The Dehn Sommerville equations hold if and only if the following equations hold:
\begin{equation}\label{eq.final}\sum_{j=k}^d(-1)^j\binom{j+1}{k+1}f_j=(-1)^df_k,\qquad k=-1,\ldots,d. \end{equation}

\end{theo}

\begin{proo}  Let $\bm f,\bm h$ be the bi-infinite vectors containing the coefficients of the f-polynomial and the h-polynomial, respectively. Following the notation shown in this paper, Equation \eqref{eq.hf} is equivalent to $\bm R_{x^{d+1}, \frac{1}{x}} \bm h = \bm h$ and the Dehn-Sommerville Equations \eqref{eq.palyndromic} are just $\bm h = \bm R_{(1-x)^{d+1}, \frac{x}{1-x}}\bm f$.

Starting with
$$\bm R_{x^{d+1}, \frac{1}{x}} \bm h = \bm h,$$

\noindent and substituting $\bm h = \bm R_{(1-x)^{d+1}, \frac{x}{1-x}}\bm f$, we get
\[
\bm R_{x^{d+1}, \frac{1}{x}} \bm R_{(1-x)^{d+1}, \frac{x}{1-x}} \bm f = \bm R_{(1-x)^{d+1}, \frac{x}{1-x}} \bm f.
\]

\noindent Now we use the multiplication rule for the product:
\[
\bm R_{(x-1)^{d+1}, \frac{1}{x-1}} \bm f = \bm R_{(1-x)^{d+1} ,\frac{x}{1-x}} \bm f.
\]

\noindent Left multiplying both sides of the equation by the inverse of $\bm R_{(1-x)^{d+1} ,\frac{x}{1-x}}$, we obtain
\[
\bm R_{(1+x)^{d+1},  \frac{x}{1+x}} \bm R_{(x-1)^{d+1}, \frac{1}{x-1}} \bm f = \bm f.
\]

\noindent Using one more time the multiplication rule, we reach:
\[
\bm R_{(-1)^{d+1}, - (1+x)} \bm f = \bm f
\]

The expression in Equation \eqref{eq.final} follows from expanding the previous matrix:
\[
\begin{bmatrix}
\ddots& \vdots & \vdots & \vdots & \vdots \\
 \hdots & \fbox{$\pm 1$} & \mp 1 & \pm 1 & \mp 1 & \hdots \\
\hdots & & \mp 1 & \pm 2 & \mp 2 & \hdots \\
\hdots & & &\pm 1 & \mp 3& \hdots \\
\hdots & & & & \mp 1 & \hdots \\
& & & & \vdots& \ddots
\end{bmatrix}
\begin{bmatrix}
\vdots \\ 0 \\ \fbox{$f_{-1}$} \\
f_0 \\
f_1 \\
f_2\\  \vdots 
\end{bmatrix}
=
\begin{bmatrix}
\vdots \\
0\\ 
\fbox{$f_{-1}$} \\
f_0 \\
f_1 \\
f_2 \\
\vdots
\end{bmatrix}.
\]





\end{proo}

\section{Some perspectives} \label{sect.pros}

As the reader may have imagined, for the (classical) Riordan group, we can also find a 2FTRM in the literature, concerning the existence of a power series whose coefficients form the so called \emph{A-sequence}. This has been deliberately omitted in the paper for the sake of brevity, but these ideas should be extended to the generalized Riordan group.

At the beginning of Section \ref{sect.comp} we have included an example, for the case $\mathbb K=\mathbb C$, for which two different formal Laurent series correspond to different expansions (at 0 and at $\infty$) of the same meromorphic function. Finding relations between one expansion and the other is a classical problem in Complex Analysis. We have intentionally abandoned the study of convergence in this setting, but maybe the framework provided here may help to approach these kinds of problem.

Also, for the (classical) Riordan group, the study of eigenproblems has received some recent attention (see, for example, \cite{CCP, LMP1, LMP2, LMP3} and the references therein). In the setting of generalized Riordan matrices that we have just described, these problems seem to be more challenging, but equally interesting (see, for example, the previous proof).

In relation to these eigenproblems, several papers in the literature address involutions in the (classical) Riordan group (see, for example, \cite{C, LMP1, LMP2} and the references therein).  We have mentioned here that the matrix $\bm J$ is an involution. It could be interesting to study, in general, the problem of finding and describing all the involutions in this new generalized context.

At the end of Section \ref{sect.riordan}, we have already pointed out that the generalized lower Riordan group appeared already in \cite{LMMPS}. In that paper, some concepts with combinatorial meaning, such as complementary and dual matrices, are discussed. These questions would also be interesting in this new framework.

\end{document}